\documentclass{llncs}
\usepackage{psfig,amssymb}

\begin{document}

\newcommand{\beq}{\begin{equation}}
\newcommand{\eeq}{\end{equation}}
\newcommand{\bea}{\begin{eqnarray*}}
\newcommand{\eea}{\end{eqnarray*}}
\newcommand{\eps}{\epsilon}
\newcommand{\Z}{{\Bbb Z}}
\newcommand{\back}{\backslash}
\newcommand{\cho}{\left(\!\begin{array}{c}}
\newcommand{\ose}{\end{array}\!\right)}
\newcommand{\mat}{\left(\!\begin{array}{cc}}
\newcommand{\rix}{\end{array}\!\right)}
\newcommand{\DET}{{\bf DET}}
\newcommand{\PP}{{\bf P}}
\newcommand{\NP}{{\bf NP}}
\newcommand{\sP}{{\bf \#P}}
\newcommand{\oa}{\overline{a}}
\newcommand{\ob}{\overline{b}}
\newcommand{\ox}{\overline{x}}

\title{Hard Tiling Problems with Simple Tiles}
\author{Cristopher Moore \inst{1,2,3} and
John Michael Robson \inst{4}}
\institute{
Computer Science Department, 
University of New Mexico, Albuquerque NM 87131
\and
Department of Physics and Astronomy, 
University of New Mexico, Albuquerque NM 87131
\and
Santa Fe Institute, 1399 Hyde Park Road,
Santa Fe, New Mexico 87501\\{\tt moore@santafe.edu}\\
\and LaBRI, Universit\'e Bordeaux 1, 351, Cours de la 
Lib\'eration, 33405 Talence, France\\
{\tt robson@labri.u-bordeaux.fr}}
\maketitle

\begin{abstract}  
It is well-known that the question of whether a given finite region
can be tiled with a given set of tiles is $\NP$-complete.  We show
that the same is true for the right tromino and square tetromino on
the square lattice, or for the right tromino alone.  In the process,
we show that Monotone 1-in-3 Satisfiability is $\NP$-complete for
planar cubic graphs.  In higher dimensions, we show $\NP$-completeness
for the domino and straight tromino for general regions on the cubic
lattice, and for simply-connected regions on the four-dimensional
hypercubic lattice.
\end{abstract}

\section{Introduction}

Tilings of the plane have long intrigued statistical physicists,
computer scientists, and recreational mathematicians.  If we have as
many copies as we like of a finite set of shapes, can we fill a given
region with them?  Tiling problems are potentially very hard: in fact,
telling whether a finite set of tiles can tile the infinite plane is
undecidable, since the non-existence of a tiling is equivalent to the
Halting Problem \cite{berger,robinson}.  For finite regions, this
problem becomes $\NP$-complete \cite{lewis}, meaning that it is just
as hard as combinatorial search problems like Traveling Salesman or
Hamiltonian Path \cite{garey}.

In fact, tiling problems can be difficult even for small sets of very
simple tiles, such as polyominoes \cite{golomb} with as few as three
cells.  Specifically, Beauquier, Nivat, R\'emila and Robson
\cite{beauquier} showed that for general regions of the square lattice
the tiling problem is $\NP$-complete for horizontal dominoes and
vertical trominoes, and for any other pair of bars where either has
length greater than 2.  This includes the straight tromino alone when
rotations are allowed.  For simply-connected regions in the plane, on
the other hand, Kenyon and Kenyon \cite{kenyon} provided a linear-time
algorithm for this problem.  When rotations are allowed, R\'emila
provided a polynomial-time algorithm for bars of length 2 and 3 in two
dimensions, both in simply-connected regions \cite{remila} and in
general \cite{remila2}.

In this paper, we give some additional examples of this kind.  For the
right tromino and square tetromino, we show that tiling a region on
the square lattice is $\NP$-complete, and that the problem of counting
how many tilings exist is $\sP$-complete, making it as hard as
enumeration problems like calculating the permanent of a matrix
\cite{papa}.  For the right tromino alone, we show $\NP$-completeness
using a more complicated construction based on 1-in-3 Satisfiability
for planar cubic graphs, but unfortunately we have no result on the
counting problem.

In higher dimensions, we consider the domino and straight tromino with
rotations allowed.  We show that this tiling problem is $\NP$-complete
for general regions in the cubic lattice, and for simply-connected
regions in the four-dimensional hypercubic lattice.

We end with a discussion of some open problems, and a discussion of
whether these results could be tightened further.

\section{Building circuits with tiles}
\label{sec:circuits}

\begin{figure}
\centerline{\psfig{file=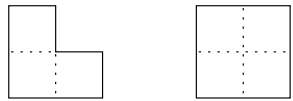,width=1.5in}}
\caption{The right tromino and square tetromino.}
\label{trisqtiles}
\end{figure}

We call the question of whether a given set of shapes can tile a given
finite subset of the square lattice, and how many such tilings exist,
the {\em existence} and {\em counting} problems.  These problems are
trivially in $\NP$ and $\sP$ respectively, since we can confirm that a
proposed tiling works in polynomial time.  In this section we will
consider the set of tiles shown in figure~\ref{trisqtiles}, namely the
right tromino and the square tetromino, with rotations allowed.

While the first hardness results about tiling relied on simulating
steps of a Turing machine from row to row, here we will instead
simulate Boolean circuits, where `wires' with two possible tilings
carry truth values, and junctions in these wires simulate logical
gates.  Similar approaches are taken in
\cite{beauquier,remila2,goles1,goles2,tdl}.  The question of whether a tiling
exists then corresponds to the canonical $\NP$-complete problem
Satisfiability, which asks whether a set of truth values for the
inputs exists that makes the output true \cite{garey}.

\begin{figure}
\centerline{\psfig{file=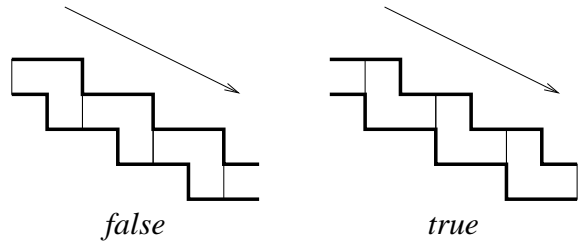,width=3in}}
\caption{The wire for tilings by the right tromino and square
tetromino.  Each wire is given an orientation, and we define its value
as true or false depending on whether the tromino occupying the
`middle' cell of each step comes from upstream or downstream.}
\label{trisqwire}
\end{figure}

Our wire is shown in figure~\ref{trisqwire}.  It moves in knight's
moves in eight possible directions across the lattice.  Given an
orientation of the wire leading from input to output, the tilings
where cells in between jogs are occupied by trominoes upstream and
downstream from them will be considered `true' and `false'
respectively.

\begin{figure}
\centerline{\psfig{file=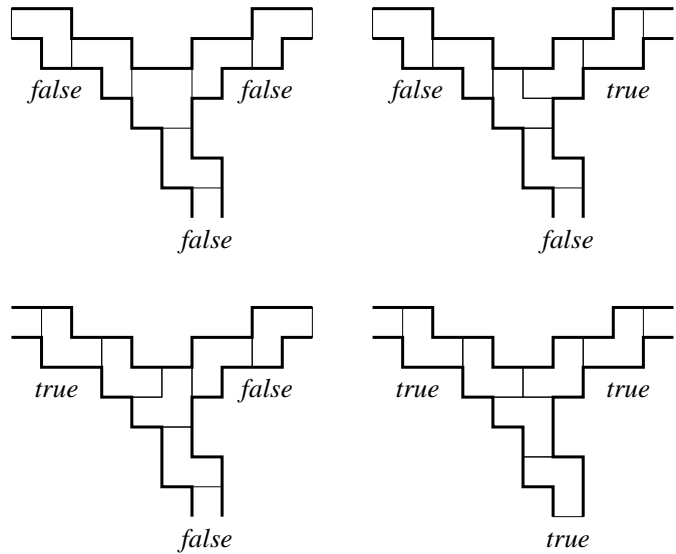,width=3.5in}}
\caption{Our AND gate.  Two input wires enter a $2 \times 2$ square
from above, and the output exits from below.  For each set of truth
values for the inputs, there is exactly one tiling as shown.}
\label{trisqand}
\end{figure}

\begin{figure}
\centerline{\psfig{file=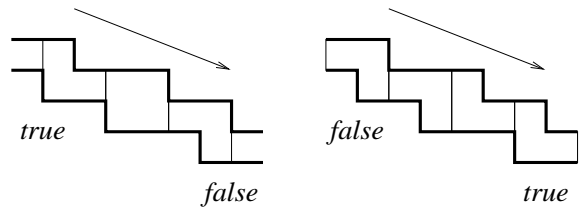,width=3in}}
\caption{Our NOT gate, which negates the signal along a wire.}
\label{trisqneg}
\end{figure}

\begin{figure}
\centerline{\psfig{file=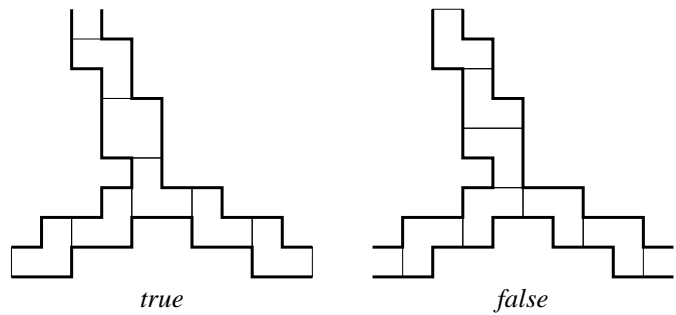,width=3.5in}}
\caption{How to split a wire (entering from above) into two copies
with the same truth value.}
\label{trisqsplit}
\end{figure}

\begin{figure}
\centerline{\psfig{file=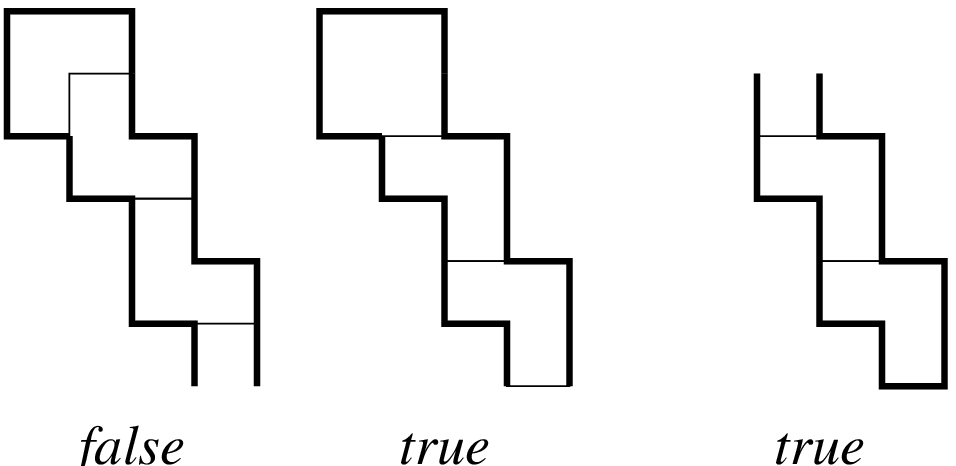,width=2.5in}}
\caption{Our variables.  The bulb at the top can be tiled in two ways,
which produce a true or false value on the wire it's connected to.  On
the right, how to end a wire so that it must be true for a tiling to
exist.}
\label{trisqvar}
\end{figure}

Given this encoding of truth values in wires, figure~\ref{trisqand}
shows an AND gate.  Each pair of truth values for the inputs $a$ and
$b$ has one and only one tiling, which causes the output to have the
value $a \wedge b$.  Figure~\ref{trisqneg} shows a NOT gate, where the
output has the opposite truth value from the input.  By combining AND
and NOT gates we can generate any Boolean function, including a
`crossover gadget' that allows us to cross wires in the plane
\cite{goldschlager}.  Since we may want to use the output of one gate
as the input of several others, figure~\ref{trisqsplit} shows how to
split a wire into two copies with the same truth value.  To create our
inputs and outputs, we need to be able to start variable wires with
either truth value, and end wires in a way that requires them to be
true.  Figure~\ref{trisqvar} shows how to do both of these.

\begin{figure}
\centerline{\psfig{file=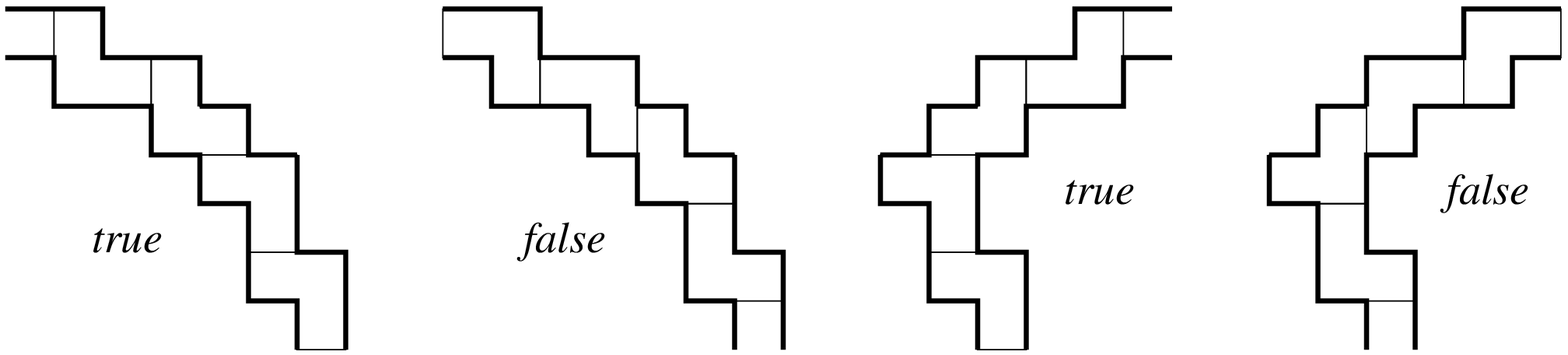,width=5in}}
\caption{How to bend wires.  On the left the wire's direction is
reflected around a diagonal line, and on the right it is rotated
$90^\circ$.  These allow us to change from any of the eight directions
to any other.}
\label{trisqturn}
\end{figure}

\begin{figure}
\centerline{\psfig{file=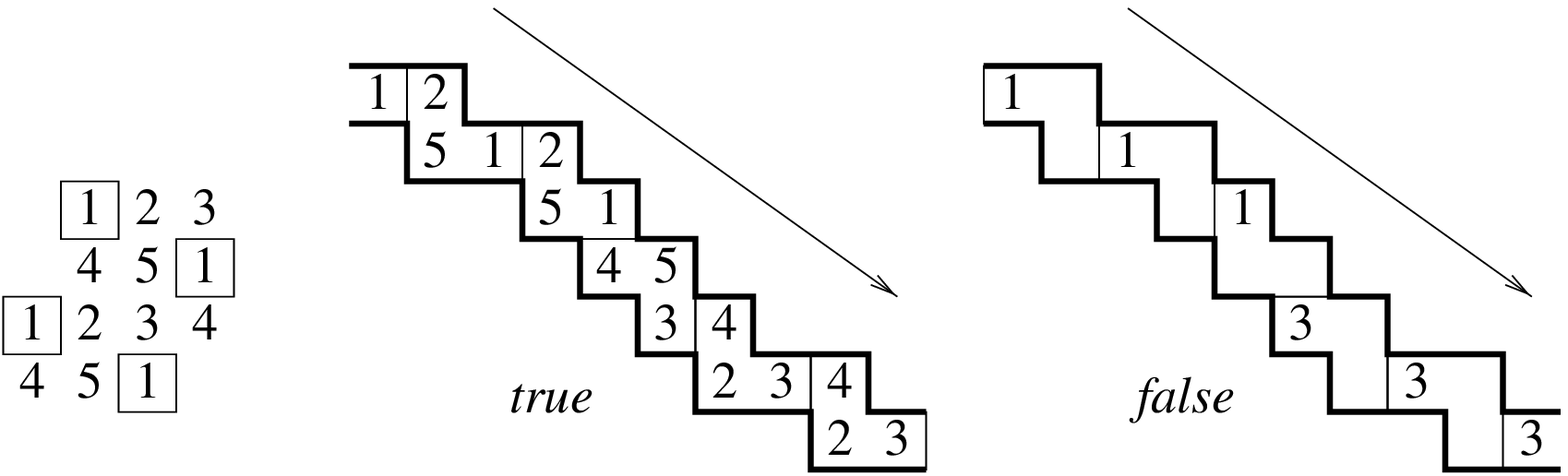,width=4.5in}}
\caption{Changing the phase of a wire.  The zig-zag shown keeps the
truth value the same, but cycles among the five colors of the
lattice.}
\label{trisqcolor}
\end{figure}

There are two more subtleties that must be dealt with.  First, it must
be possible to bend wires among their eight possible orientations.  We
can do this by rotating and reflecting a wire's direction as in
figure~\ref{trisqturn}; these turns are derived simply by setting one
of an AND gate's inputs to be true.  Secondly, a wire of a given
orientation has a knight's move periodicity, so that if we color the
cells of the lattice as in figure~\ref{trisqcolor} each step can only
land on one of the five colors.  This gives each wire a kind of phase,
and the zig-zag shown in the figure allows us to change this phase so
that a wire's truth value can be delivered to cells of any color.
Then the output of any gate can become the input of any other, as long
as the gates are separated widely enough.  It is easy to show that
this separation need only be polynomially large.

It is clear from this construction that for any Boolean circuit, we
can construct a region of the plane such that the right tromino and
square tetromino can tile it if and only if the circuit is
satisfiable.  Moreover, this reduction is {\em parsimonious}, i.e.\
the mapping between solutions of the two problems is one-to-one, so
the number of solutions is the same in both cases.  Thus the question
of whether any tilings exist, and how many there are, are equivalent
to Satisfiability and \#Satisfiability, which are $\NP$-complete and
$\sP$-complete respectively \cite{garey,papa}.

\section{Right Trominoes Alone}

We now improve on the result of the previous section, and show that
the existence problem is $\NP$-complete for the set of tiles
consisting only of the right tromino (with rotations allowed).  The
main reason that this is more difficult than for the tromino and
square tetromino is that every tile covers three squares.  If we use
the same wires as before, where the truth value depends on the
protrusion, or lack of it, of a single site, then we can only
represent gates where the output is a linear function of the inputs
mod 3.  While this prevents us from coding Boolean circuits directly,
we can still represent a version of Satisfiability which we will show
is $\NP$-complete.

\begin{figure}
\centerline{\psfig{file=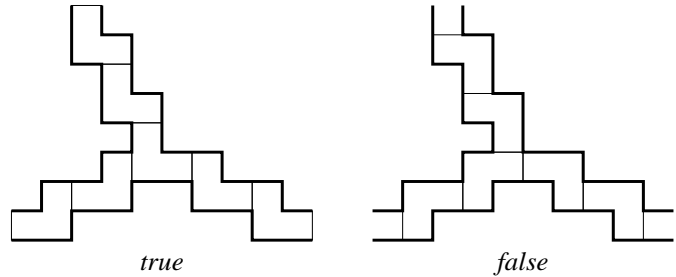,width=3.5in}}
\caption{The variable nodes for our reduction of cubic planar 1-in-3
Satisfiability to tromino tiling.  All three outgoing wires have the
same truth value.  We can also use this node to enforce equality
between three incoming wires.}
\label{trovar}
\end{figure}

First, we represent variables with nodes as in Figure~\ref{trovar}, in
which three outgoing wires are constrained to be all true or all
false.  (The same node can be used to enforce equality among three
incoming wires.)  This gate is essentially the wire splitter of
Figure~\ref{trisqsplit} with the NOT gate removed from the top wire.

\begin{figure}
\centerline{\psfig{file=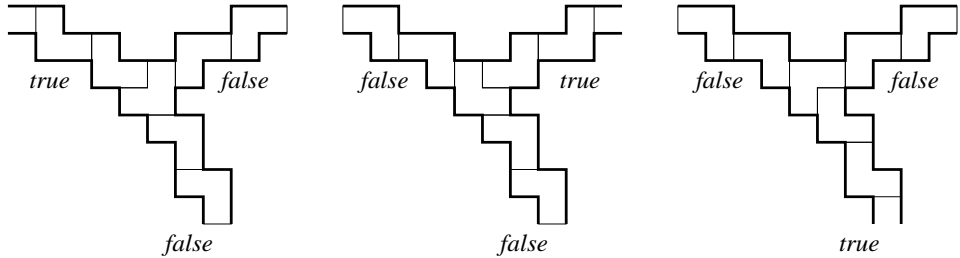,width=5in}}
\caption{The clause nodes for our reduction of cubic planar 1-in-3
Satisfiability to tromino tiling.  Exactly one of the three incoming
wires must be true.}
\label{tro1in3}
\end{figure}

Secondly, we represent clauses with the nodes of Figure~\ref{tro1in3},
which require that exactly one of three incoming wires be true, and
that the other two be false.  This is a slightly modified version of
the AND gate of Figure~\ref{trisqand}.  If exactly one wire is true,
the tilings of the wires cover exactly one square of the central $2
\times 2$ square, leaving a space for one right tromino.  Otherwise
the number of squares left in the central square is not a multiple of
3 and so cannot be tiled by trominoes.

If we consider a bipartite graph of these two types of nodes
connected with wires as in the previous section, then the tiling
problem becomes a variant of 1-in-3 Satisfiability.  That is, we have
a set of variables and a set of three-variable clauses, and we want to
know if there is an assignment of truth values to the variables such
that exactly one of the three variables in each clause is true.

However, here we have a restricted case of this problem in which the
expression's graph is planar, there is no negation, and each variable
occurs in exactly three clauses since our variable nodes have three
wires emanating from them.  Thus, to complete our proof we need to
show that Cubic Planar Monotone 1-in-3 Satisfiability is
$\NP$-complete, and we do this in the next section.

\subsection{Cubic Planar Monotone 1-in-3 Satisfiability}
\label{cubic}

LaRoche \cite{laroche} showed that Planar Monotone 1-in-3 SAT is
$\NP$-complete.  This is the restriction of Monotone 1-in-3 SAT to
expressions whose graph $G$ is planar, where $G$ is defined with a
vertex for each variable and each clause of the expression, and an
edge joining each clause to each variable appearing in it.

We will show that this problem remains NP-complete when restricted to
expressions where every variable has exactly three occurrences, i.e.\
those whose graph is cubic.  (An alternate formulation is to consider
a hypergraph where each clause is a node and each variable is a
hyperedge joining three nodes.  Since each true variable ``covers''
three clauses and every clause must be covered exactly once, we can
then regard this problem as Exact Cover by 3-Sets in the case where
the associated hypergraph is planar.)

While the reduction from this version of SAT to the tiling problem we
gave in the previous section is parsimonious, this does not show that
the counting problem is $\sP$-complete, since we do not have such a
result for Cubic Planar Monotone 1-in-3 SAT.  We leave this as an open
problem for the reader.

Our reduction is from unrestricted Planar Monotone 1-in-3 SAT.  In the
diagrams of this section we use the same convention as in
\cite{laroche} that circles represent variables and squares represent
clauses.

The reduction proceeds in two stages. The first stage produces an
expression $E'$ which is planar and satisfiable if and only if the
original was, and where each variable has three or fewer occurrences.
The second stage, rather more intricate, adds new clauses and
auxiliary variables producing an expression $E''$ where all variables
have exactly three occurrences.

\subsection{The reduction}
\label{sec:reduction}

\begin{figure}
\centerline{\psfig{file=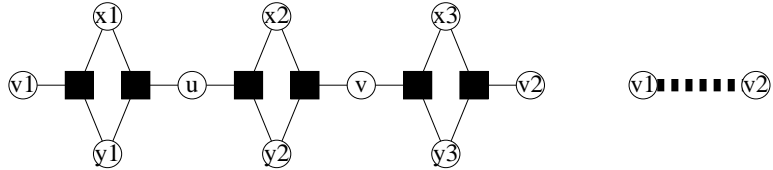,width=4in}}
\caption{A chain that forces two variables $v_1$ and $v_2$ to be equal, 
and our abbreviation for it.}
\label{chain}
\end{figure}

We start by introducing some simple components to be used in the first
stage.  An {\em equality verifier} for two variables $v_1$ and $v_2$
consists of two new variables $x$ and $y$ together with the two
clauses $(v_1,x,y)$ and $(v_2,x,y)$.  Clearly this can be satisfied if
and only if $v_1=v_2$; moreover, if $v_1$ and $v_2$ are false, then
exactly one of $x$ and $y$ is true.  A {\em chain} between two
variables $v_1$ and $v_2$ consists of new variables $u$ and $w$ and
three equality verifiers: one between $v_1$ and $u$, one between $u$
and $w$, and one between $w$ and $v_2$.  Again a chain forces
$v_1=u=v=v_2$ in any satisfying assignment and when $v_1$ and $v_2$
are false, the choice of $x$ or $y$ being true is independent for the
three equality verifiers.  A chain is shown in Figure~\ref{chain}.

\begin{figure}
\centerline{\psfig{file=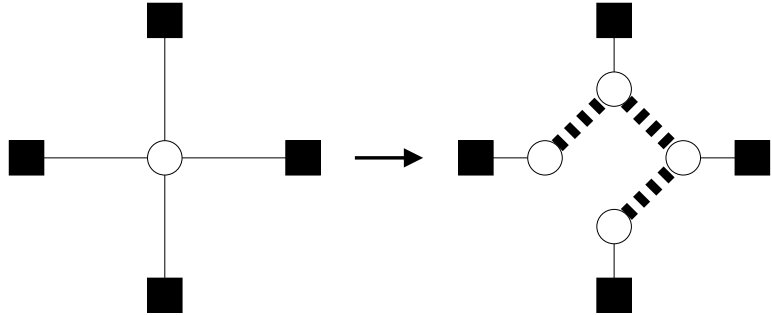,width=4in}}
\caption{Replacing variables with $k > 1$ occurrences with $k$ variables
connected by chains.  One link is omitted arbitrarily, say between
$v_{i,1}$ and $v_{i,k}$, so as not to form a complete cycle.}
\label{replace}
\end{figure}

Given an expression $E$ which is an instance of Planar Monotone 1-in-3
SAT, we choose a planar embedding of its graph $G$.  For each variable
$v_i$ that occurs in $k > 1$ clauses, we replace it with a set of
variables $v_{ij}$ for $j=1 \ldots k$, one for each occurrence.  We
then link these $k$ variables together with chains as in
Figure~\ref{replace}.  We call the resulting expression $E'$ and its
planar graph $G'$.  Note that this sequence of links between the
$v_{ij}$ is not completed to form a cycle; we omit an arbitrarily
chosen link, say between $v_{i,1}$ and $v_{i,k}$.  Therefore, any
boundary between two faces that includes a chain in fact includes
chains from two different variables that appear in the same clause, a
fact which we will use below.

This completes the first stage of the reduction.  Clearly $E'$ is
satisfiable if and only if $E$ was, and every variable has three or
fewer occurrences.  The variables with less than three occurrences are
the $x$ and $y$ variables of each equality verifier, the intermediate
$u$ and $w$ variables of each chain, the two unlinked variables
$v_{i,1}$ and $v_{i,k}$ for each $v_i$ with $k > 1$ occurrences, and
the original variables $v_i$ which have only one occurrence in $E$.
The total number of missing occurrences is a multiple of 3 since each
variable must have three occurrences and each clause absorbs three.

We now consider the faces of $G'$.  Our goal is to associate each
missing occurrence to a face so that each face has a multiple of three
additional occurrences along its boundary.  There are some faces
bordered by the $v_1,x,v_2,y$ of an equality verifier; we will call
these faces trivial, and assign no missing occurrences to them.  
The remaining assignments are made as follows.

Since the $x$ and $y$ variables of the chains have only one face
adjacent to them, we start with each non-trivial face having a certain
number of additional occurrences mod 3.  We now note that every
non-trivial face has at least one chain along its boundary.  Each
chain has two variables $u$ and $v$ with one occurrence each that can
be assigned either to this face or to one adjacent to it.  In this
way, we can transfer 0, 1 or 2 (mod 3) additional occurrences from
this face to the next one.  If we define a tree of faces (i.e.\ a
spanning tree of the dual graph) we can start at the leaves and change
the number of additional occurrences in each face to $0 \bmod 3$,
until we reach the root, where the number of occurrences left is also
$0 \bmod 3$ since the total is a multiple of 3.

We use a similar strategy to control, for each face, how many of the
additional occurrences entering it are true, assuming that $E'$ is
satisfiable.  Recall that along the boundary between any pair of
adjacent non-trivial faces are chains corresponding to two different
variables that appear in the same clause.  In any satisfying
assignment, at least one of these is false, and so in each link of its
chain we can take either $x$ or $y$ to be true.  In this way, we can
transfer 0, 1 or 2 (mod 3) true variables from one face to another.
If we again define a tree of faces, we can start at the leaves and
change the number of true additional occurrences in each face to $0
\bmod 3$, until we reach the root.  At the root face, the number of
true additional occurrences is simply the total mod 3.  This is
$-(\mbox{the number of clauses of }E') \bmod 3$, since variables
provide true values in multiples of 3 and each clause of $E'$ absorbs
one of these.

We will now construct our final expression $E''$, and its graph $G''$,
by gathering the additional occurrences around each face into a {\em
gadget} as described in Section~\ref{sec:gadgets}.  This gadget will
enforce that the correct number mod 3 of the additional occurrences
are true.  In what follows, we will sometimes write 0 and 1 for true
and false respectively.

\subsection{Gadgets}
\label{sec:gadgets}

The purpose of a gadget is to be satisfiable provided the correct
number mod 3, say $c$, of its $3m$ input edges are true.  Thus we will
complete the construction of $G''$ by placing a gadget in each face of
$G''$.

A gadget is composed of a sequence of {\em optional switches} that
allow the true and false values to be sorted in any order, and checked
in groups of three.  Specifically, we wish for the rightmost $c$
groups of three variables to have one true input each, and for all
other groups of three variables to be either all true or all false.
To do this, we terminate the rightmost $c$ groups with a single clause
vertex, and all other groups with a {\em 3-way equality verifier} as
defined below.

An optional switch has two entry edges and two exit edges.  It can be
satisfied in all cases where the two values on the exit edges are the
same as the two values on the entry edges in either order.  It uses
two sub-gadgets, the {\em partial switch} which behaves the same way
except that it cannot be satisfied if both entry edges are true, and
the {\em 3-way equality verifier}, which requires all its inputs to
have the same truth value.

The partial switch itself has a sub-gadget which we call a {\em
triangle}.  This is a component with three variable vertices
$(a,b,c)$, each with one edge to the exterior, which requires that
exactly one of $a$, $b$ and $c$ is true.  In other words, it acts just
like a clause, except that it uses two of each variable's occurrences.
Note that we cannot achieve this simply by doubling the clause
$(a,b,c)$ since this would violate planarity.

\begin{figure}
\centerline{\psfig{file=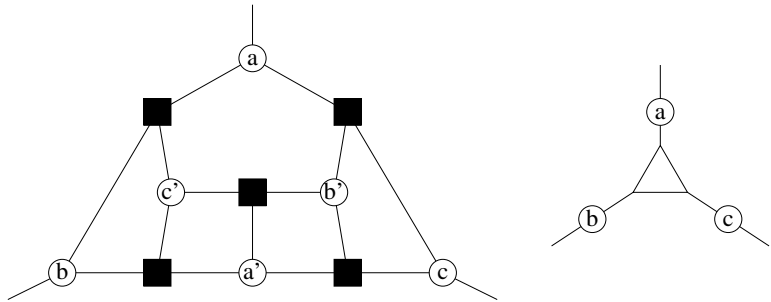,width=4in}}
\caption{A triangle and our abbreviation for it.}
\label{triangle}
\end{figure}

\begin{figure}
\centerline{\psfig{file=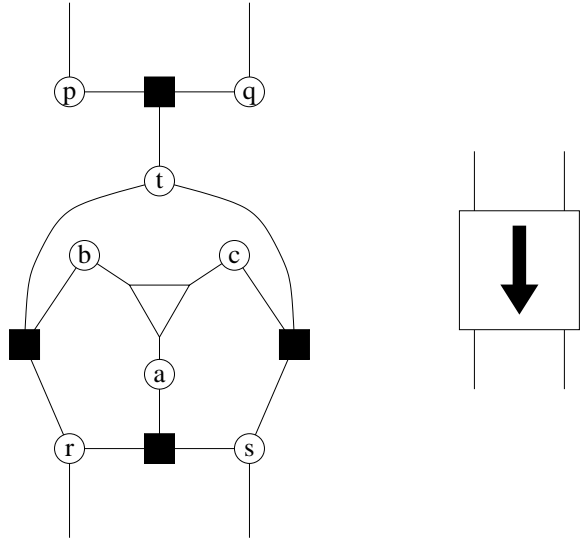,width=3in}}
\caption{A partial switch and its abbreviation.}
\label{partial}
\end{figure}

The triangle has three internal variables $a'$, $b'$, and $c'$, and
the clauses $(a,b,c')$, $(b,a',c')$, $(c,a',b')$, $(a,c,b')$ and
$(a',b',c')$.  Clearly this is satisfied by taking $a'=a$, $b'=b$ and
$c'=c$, and the central clause then ensures that exactly one of these
is true.  A triangle is shown in Figure~\ref{triangle} together with
the symbol used to denote it in larger components.

The partial switch is shown in Figure~\ref{partial}.  It has two input
variables $p,q$ and two output variables $r,s$.  It contains internal
variables $a, b, c$ and $t$, a triangle $(a,b,c)$, and the clauses
$(p,q,t)$, $(b,r,t)$, $(c,s,t)$, $(a,r,s)$.  If $p=q=0$, the partial
switch is satisfied by taking $a=t=1$ and $b=c=r=s=0$.  If $p \ne q$,
we can take either $r=c=1$ and $a=b=s=t=0$ or $s=b=1$ and $a=c=r=t=0$,
and in either case $r \ne s$.  Thus the partial switch copies the two
inputs to the two outputs in either order, except when both inputs are
true.

\begin{figure}
\centerline{\psfig{file=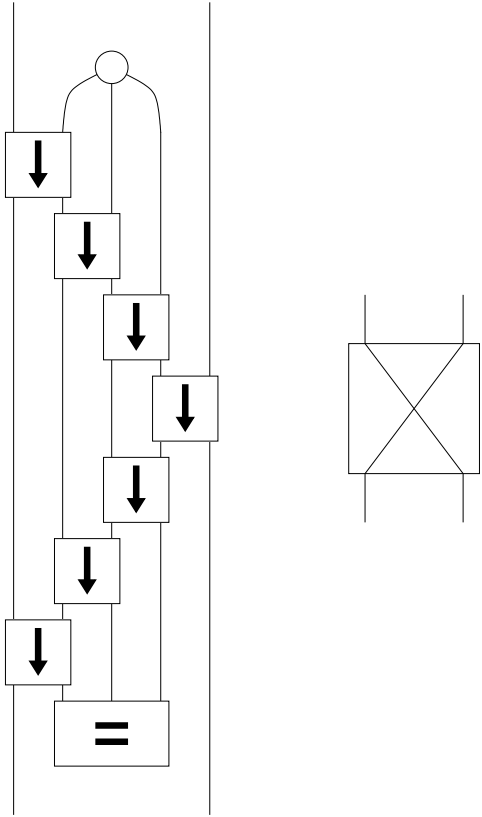,width=2.5in}}
\caption{An optional switch and its abbreviation.}
\label{switch}
\end{figure}

Note that each partial switch uses two occurrences each of its input
variables, and only one each of its output variables.  Thus partial
switches can be strung together in series.  We use this to build the
full optional switch shown in Figure~\ref{switch}.  This is composed
of a new variable vertex (which may always be false) placed between
the two entry edges followed by two sequences of partial switches.
These permit a single true value to be transmitted from entry 1 to
exit 2 or vice versa.  This leaves three edges between the two exit
edges and these three edges are then passed into a 3-way equality
verifier.

\begin{figure}
\centerline{\psfig{file=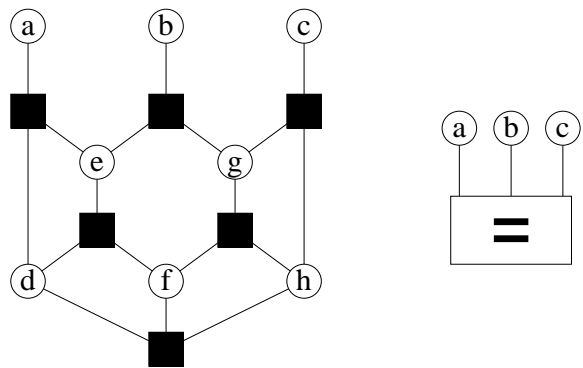,width=3in}}
\caption{A 3-way equality verifier and its abbreviation.}
\label{equality}
\end{figure}

\begin{figure}
\centerline{\psfig{file=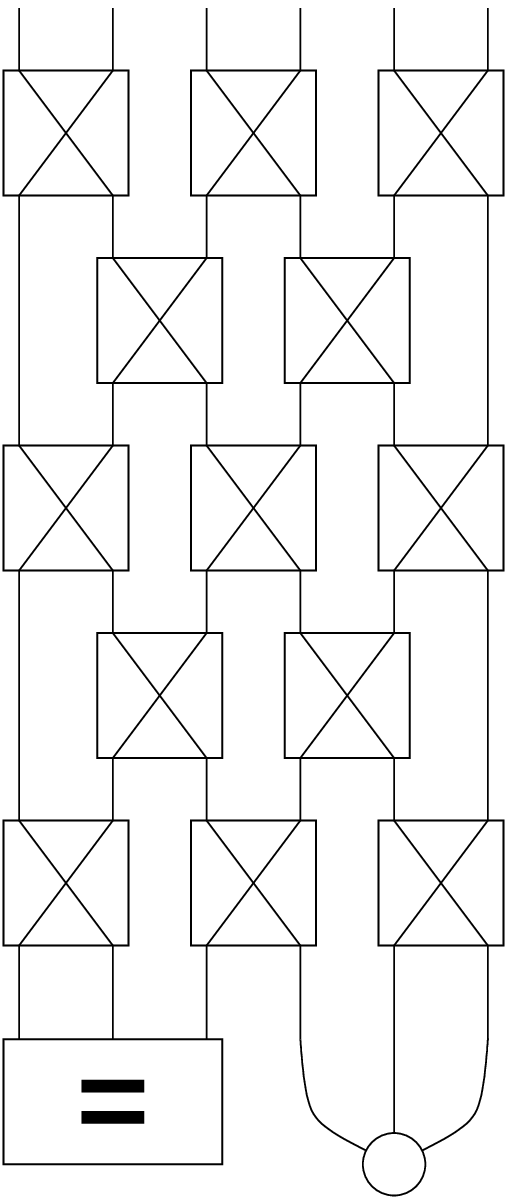,width=1.5in}}
\caption{A gadget which allows 6 inputs to be sorted and confirms 
that $1 \bmod 3$ of them are true.}
\label{gadget}
\end{figure}

For a 3-way equality verifier in the tromino tiling problem, we can
simply use the variable node of Figure~\ref{trovar} with three
incoming wires rather than three outgoing ones.  However, to complete
our proof for Cubic Planar Monotone 1-in-3 SAT, we construct a
verifier as in Figure~\ref{equality}.  It has three input edges from
variables $a,b,c$ and five internal variables $d,e,f,g,h$.  It
contains the clauses $(a,d,e)$, $(b,e,g)$, $(c,g,h)$, $(d,e,f)$,
$(f,g,h)$, and $(d,f,h)$.  If $a=b=c=0$, we take either $d=g=1$ and
$e=f=h=0$ or $e=h=1$ and $d=f=g=0$.  If $a=b=c=1$ we take $f=1$ and
$d=e=g=h=0$.

Finally, to construct a gadget for $3m$ inputs of which $c \bmod 3$
are to be true, we form a sorting network from a polynomial number of
optional switches that allows the inputs to be sorted in any order.
We then divide the outputs of the sorting network into triples, and
add a suitable checker to each triple.  This checker will be a single
clause for $c$ of the triples and a 3-way equality verifier for the
rest.  An example is shown in Figure~\ref{gadget}, where $m=2$ and
$c=1$.

After the appropriate gadget is added to each face to absorb its
additional occurrences, it is clear that $E''$ is satisfiable if and
only if $E$ is.  This completes our proof that Cubic Planar Monotone
1-in-3 Satisfiability, and therefore the existence problem for tilings
by the right tromino, is $\NP$-complete.

\section{Monotone circuits in three and four dimensions}

\begin{figure}
\centerline{\psfig{file=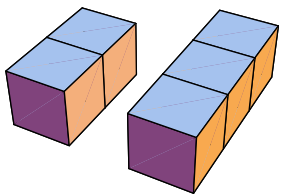,width=1.5in}}
\caption{The domino and straight tromino in three dimensions.}
\label{domtri}
\end{figure}

In this section, we construct wires and gates in three dimensions with
the domino and straight tromino, or (if you prefer) the dicube and
straight tricube, shown in figure~\ref{domtri}.  The construction is
similar to that in Section~\ref{sec:circuits}, except that it will be
limited to {\em monotone} Boolean circuits, whose outputs are
non-decreasing functions of their inputs.

\begin{figure}
\centerline{\psfig{file=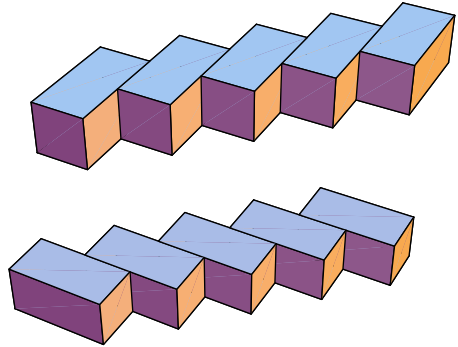,width=2.5in}}
\caption{The wire for tilings by the domino and straight tromino.  It
can be tiled by dominoes in either of two ways, such that the
downstream end of each domino lies on either an odd cell (true) or an
even cell (false).}
\label{domtriwire}
\end{figure}

\begin{figure}
\centerline{\psfig{file=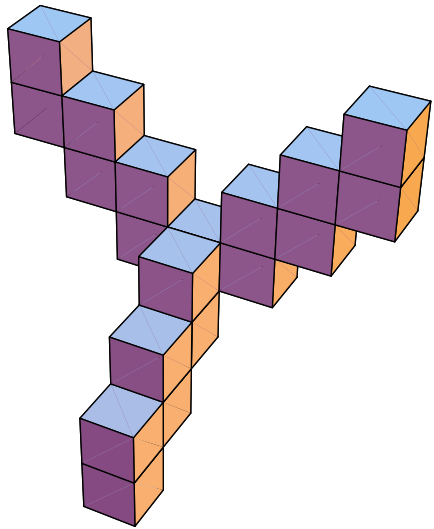,width=2.5in}}
\caption{This configuration serves as an AND or OR gate, depending on
whether the central cell is even or odd.}
\label{domtrigate}
\end{figure}

Our wires will consist of zig-zags as in figure~\ref{domtriwire},
which can be tiled by dominoes in either of two ways.  As before, we
give each wire an orientation.  Then if we two-color the lattice as a
three-dimensional checkerboard, dividing it into odd and even cells,
we define true and false wires as those where the downstream end of
each domino lies on an odd or even cell respectively.  Note that this
allows wires to turn quite freely.  Then the reader can easily check
that the configuration in figure~\ref{domtrigate}, with two inputs
entering from above and an output exiting from below, acts as an AND
gate if its central cell is even, and an OR gate if its central cell
is odd.  There is a unique tiling for each truth value of the inputs;
in one of these four tilings the central crossbar is filled with a
tromino, while the other tilings consist entirely of dominoes.

With AND and OR gates, we can build any monotone Boolean circuit.  It
seems to be difficult to build a NOT gate with these two tiles, and we
conjecture that this is impossible.  However, monotone circuits are
enough for our purposes, since we can use De Morgan's laws
$\overline{a \wedge b} = \oa \vee \ob$ and $\overline{a \vee b} = \oa
\wedge \ob$ to move negations up through our gates
\cite{goldschlager}, and thus convert any Boolean circuit with
variables $x_1, \ldots, x_n$ to a monotone one with variables $x_1,
\ox_1, \ldots, x_n, \ox_n$.

To make a one-to-one correspondence between satisfying assignments of
this circuit and the original one, we just have to add the
non-monotone condition that exactly one of $x_i$ and $\ox_i$ be true
for each $i$.  We can do this by connecting these two variables with a
zig-zag wire; using our definition of orientation, this wire's tiling
will come out as true at one end and false at the other.  Finally, we
require the output wire to be true by ending it on an odd cell.

\begin{figure}
\centerline{\psfig{file=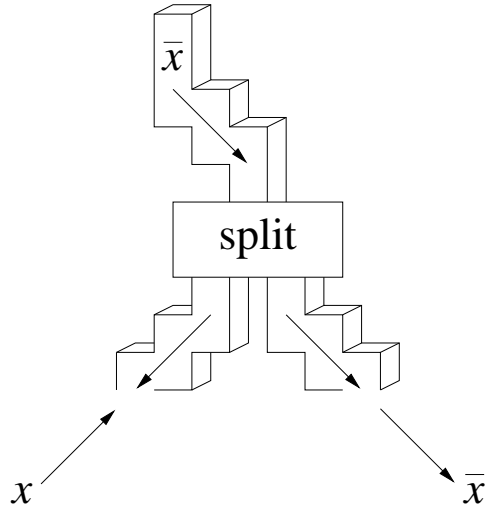,width=2.5in}}
\caption{If we had a wire splitter, we could make a NOT gate by
feeding a wire with value $x$ into one of its outputs in reverse, and
connecting its input into a region which can be tiled either true or
false.  Since our definition of truth value is negated if we reverse
orientation, the other output will then carry the negated value
$\ox$.}
\label{split}
\end{figure}

To complete the argument, we need to be able to split wires.  However,
if we can do this, we can construct a NOT gate as shown in
figure~\ref{split}.  Since our definition of truth is negated if we
reverse the orientation of a wire, if we attach one output of a
splitter backwards to a wire with value $x$, and attach its input to a
variable that can be true or false, the other output will yield $\ox$.
Thus if we are correct that we can't build a NOT gate, we can't build
a wire splitter either.

Since a perfect wire splitter doesn't seem to exist, we consider
instead a ``dirty splitter'' where the outputs are equal to or less
than the inputs; for instance, an OR gate in reverse has this
property.  Since the circuit is monotone, this can create false
negatives but not false positives, so this ``dirty circuit'' can have
a true output if and only if the original circuit can.  Thus we have
shown how to convert a monotone Boolean circuit to a region of the
cubic lattice that can be tiled with dominoes and straight trominoes
if and only if the circuit is satisfiable, and the existence problem
for such tilings is $\NP$-complete.

We would like to say that the counting problem for such tilings is
$\sP$-complete, but using dirty splitters means that in some cases
there will be more than one tiling for a given satisfying assignment.
Specifically, if a wire coming out of a dirty splitter is true but is
not essential to making the output true, we can tile it false and the
circuit will still be satisfied.  Thus the mapping between assignments
and tilings is not one-to-one, and our reduction from Satisfiability
to the existence problem is not parsimonious.  While it seems likely
that counting these tilings is $\sP$-complete, it is not clear how to
prove this.

It unlikely that the number of dimensions in this construction can be
reduced from three to two.  While planar circuits and monotone
circuits can both simulate Boolean circuits in general
\cite{goldschlager}, circuits which are both planar and monotone
cannot.

Finally, we note that we can convert any finite region $R$ in the
cubic lattice into a simply-connected region $R'$ in the
four-dimensional hypercubic lattice with an equivalent tiling problem
\cite{personal}.  Simply choose a set of sites $S$ such that $R \cup
S$ is simply-connected, two-color them, and add each one along with an
additional site adjacent to it in the fourth dimension, alternating up
and down according to the two-coloring.  Then these sites can only be
filled by dominoes aligned along the fourth dimension, and only the
sites in $R$ are left.  Thus the tiling problem for dominoes and
straight trominoes is also $\NP$-complete for simply-connected regions
in the four-dimensional cubic lattice.

\section{Discussion}

For the right tromino and square tetromino on the square lattice, or
the right tromino alone, we have shown that the existence problem for
tiling regions of the square lattice is $\NP$-complete.  In the former
case, we have shown that the counting problem is $\sP$-complete.  For
the domino and straight tromino, we have shown that the existence
problem is $\NP$-complete for general regions in the cubic lattice,
and for simply-connected regions in the four-dimensional hypercubic
lattice.

Intuitively, the $\NP$-completeness of tiling with polyominoes with
three or more cells comes from that fact that, while telling whether a
graph can be covered with dimers is related to the Perfect Matching
problem and can be solved in polynomial time, telling whether it can
be covered by trimers is $\NP$-complete \cite{garey}.  While there is
no strict connection between a system's computational complexity and
its statistical mechanics, we might expect that tiling problems for
tiles consisting of three or more cells are generally not exactly
solvable in two or more dimensions (although the statistics of
triangular trimers on the triangular lattice can be solved with a
Bethe ansatz \cite{trimers}).

For dominoes the existence and counting problems have an especially
elegant solution.  The number of ways to tile a region of a regular
lattice with dominoes can be calculated by expressing it as a
determinant of a modified adjacency matrix \cite{kasteleyn,propp}.
This puts the existence and counting problems for the domino in the
class $\DET \subset \PP$ and allows the statistics of the dimer model
to be solved exactly.

It would be interesting to find other small sets of tiles for which
tiling problems are hard, or find larger sets of tiles for which
existence and/or counting are in $\PP$.  As noted above, in two
dimensions the existence problem for dominoes and straight trominoes
is in $\PP$ when rotations are allowed \cite{remila,remila2}

Another solvable case is that of a single polyomino without rotation,
or one which is rotationally symmetric.  Then there is at most one
tiling of a given finite region, which we can find by scanning from
left to right and top to bottom and adding a tile to each unoccupied
cell we meet.  Beauquier and Nivat \cite{bn} showed that this case of
the tiling problem is decidable for the infinite plane as well.
Finally, we note that the existence of an isohedral tiling of the
infinite plane by a single polyomino --- one whose symmetry group acts
transitively on the set of tiles --- was recently shown to be
decidable by Keating and Vince \cite{isohedral}.

Some open questions include:
\begin{enumerate}
\item Is the counting problem for the right tromino on the square
lattice, or for the domino and straight tromino on the cubic lattice,
$\sP$-complete?
\item Is the existence problem for the right tromino in $\PP$ for
simply-connected regions in the square lattice?
\item Is the existence problem for dominoes and trominoes in $\PP$ for
simply-connected regions in the cubic lattice? 
\end{enumerate}

{\bf Acknowledgements.}  We are grateful to Richard Kenyon, Michael
Lachmann, Mark Newman, Mats Nordahl, Laurent Vuillon, and Eric
R\'emila for helpful conversations.  C.M.\ also thanks Robert Cori for
an invitation to \'Ecole Polytechnique, to Michel Morvan for an
invitation to Universit\'e Paris 7 (Jussieu) where the authors met,
and to Molly Rose and Spootie the Cat for warmth and friendship.

\end{document}